\baselineskip=14pt
\font\Bf=cmbx12
\def\epsi{\varepsilon}
\def\x{x^*}
\def\y{y^*}
\def\z{z^*}
\def\U{u^*}
\def\v{v^*}
\def\B{B^*}
\def\X{X^*}

\def\<{\langle}
\def\>{\rangle}
\def\tto{\,\scriptstyle$\lower1pt\vbox{\hbox{$\to$}\kern-10pt\hbox{$\to$}}$\,}

\centerline{\Bf Regular Maximal Monotone Multifunctions and Enlargements}

\medskip
\centerline{\bf Andrei Verona}
\centerline{Department of Mathematics, California State University,}
\centerline{Los Angeles, CA 90032, USA}
\centerline{e-mail: averona@calstatela.edu}

\medskip
\centerline{\bf Maria Elena Verona}
\centerline{Department of Mathematics, University of Southern California}
\centerline{Los Angeles, CA 90089-1113, USA}
\centerline{e-mail: verona@usc.edu}

\bigskip
\centerline{\bf Dedicated to Stephen Simons}

\bigskip\noindent
{\bf Abstract.} In this note we use recent results concerning the sum theorem for maximal monotone multifunctions
in general Banach spaces to find new characterizations and properties of regular maximal monotone multifunctions
and then use these to describe the domain of certain enlargements.



\bigskip\noindent
{\bf Introduction}

\medskip
Throughout this note $X$ will denote a Banach space with topological dual $\X$. For
$(x,\x)\in X\times \X$, $\<x,\x\>$ (or $\<\x,x\>$) will denote the natural evaluation map.
The unit ball in $\X$ is denoted $\B$. If $T:X\tto \X$ is a multifunction then
$G(T)=\{(x,\x)\in X\times \X: \x\in T(x)\}$ is called its {\it graph}, while
$D_T=\{x\in X:T(x)\ne \emptyset\}$ is called its {\it domain}. The multifunction
$T$ is called {\it monotone} if $\<\x-\y,x-y\>\ge 0$ for all $(x,\x),(y,\y)\in G(T)$;
$T$ is called {\it maximal monotone} if there exists no monotone multifunction $S:X\tto \X$ such
that $G(S)$ strictly contains $G(T)$. The theory of maximal monotone multifunctions in
reflexive spaces is now more or less complete, due to Rockafellar's results from the
1960s and 1970s. However, in general Banach spaces, the theory is much more complicated.
In order to extend Rockafellar's results to general Banach spaces, different authors
have introduced several classes of maximal monotone operators and proved some of his results in these
particular cases (see [7], [8] for a detailed description of these classes and results
and also [11], [12] and [13] for properties of regular maximal monotone multifunctions).

In this note we use recent results concerning the sum theorem for maximal monotone multifunctions
in general Banach spaces to find new characterizations for regular maximal monotone multifunctions
(Theorem 2, Corollary 3, Theorem 4). We also prove that the sum of a regular maximal monotone multifunction and
a bounded maximal monotone multifunction is maximal. In the last part we characterize a type of enlargements (Theorem 7),
give a characterization of regularity in terms of certain enlargements (Theorem 11), and show that the domain
of another type of enlargements is contained in the closure of the domain of the initial multifunction (Theorem 9).

It is worth mentioning that if a sum theorem for maximal monotone multifunction in general Banach spaces was proved then any
maximal monotone multifunction would be regular and therefore all our results in this note would be true for any maximal
monotone multifunction.

\medskip\noindent
{\bf Regular maximal monotone multifunctions}

\medskip\noindent
In [11] we introduced the following number
$$L(x,\x,T)=0\vee\sup\{{\<\x-\y,y-x\>\over \|x-y\|}:\,y\ne x, (y,\y)\in G(T)\}$$
(here $a\vee b=\max\{a,b\}$) and proved that $L(x,\x,T)\le \inf\{\|\y-\x\|: \y\in T(x)\}=d(\x,T(x))$. (When $T$ is
the subdifferential of a proper, convex, lower semicontinuous function this number was
considered by Simons [9].) The maximal monotone multifunction $T$ is called {\it regular} if
$L(x,\x,T)=\inf\{\|\y-\x\|: \y\in T(x)\}=d(\x,T(x))$ for any $(x,\x)\in X\times \X$. Here is a list
with some relevant facts about regular maximal monotone multifunctions:

\smallskip - If $T$ is a regular maximal monotone multifunction then $\overline {D_T}$ is convex [11, Theorem 2].

\smallskip - Any maximal monotone multifunction in a reflexive Banach space is regular [11, Corollary 1 (2)].

\smallskip - The subdifferential of any proper convex lower semicontinuous function is regular [9, Theorem 6].

\smallskip - If $T$ is a regular maximal monotone multifunction and $x\in {\overline D_T}$ then
$T$ is locally bounded at $x$ if and only if $x\in {\rm int} D_T$ [11, Corollary 3].

\smallskip - If $T$ is a linear (possibly discontinuous) and maximal monotone multifunction then $T$ is
regular [12, Proposition 3.2].

\smallskip - If $T$ is a strongly-representable maximal monotone multifunction then $T$ is regular (see Remark 7 in [15]).
Since any maximal monotone multifunction of type (NI) is strongly-representable (Proposition 26 in [15]) it follows that
any maximal monotone multifunction of type (NI) is regular (recently it was proved [4, Theorem 1.2] that the class of
strongly-representable monotone multifunctions coincides with the class of maximal monotone multifunction of type (NI)).
In particular any maximal monotone multifunction of type (D)
is regular (since any maximal monotone multifunction of type (D) is of type (NI), [8, Theorem 36.3(a)]).

\smallskip
We shall prove in this note that:

\smallskip - If $T$ is a maximal monotone multifunction and $D_T$ is either closed and convex or has non-empty
interior then $T$ is regular (Corollary 3).

\smallskip - If $T$ is a maximal monotone multifunction then $T$ is regular if and only if it is dually strongly
 maximal (Theorem 4).

\medskip\noindent
{\bf Theorem 1.} Let $T:X\tto \X$ be a regular maximal monotone multifunction and $S:X\tto \X$
be a bounded maximal monotone multifunction. Then $T+S$ is maximal monotone.

\medskip\noindent
{\bf Proof.} In view of the Debrunner-Flor Extension Theorem (see [5, Lemma 1.7]), $D_S=X$. Let $M$ be such that
$\|\U\|\le M$ for any $(y,\U)\in G(S)$. Let $(x,\x)\in X\times \X$ be monotonically related to $G(T+S)$, that is
$$\<\y+\U-\x,y-x\>\ge 0\ \ {\rm for\ any}\ \ (y,\y)\in G(T),\ (y,\U)\in G(S).$$
Then
$$\<\y-\x,x-y\>\le\<\U,y-x\>\le\|\U\|\,\|y-x\|\le M\|y-x\|$$
from which it follows that $L(x,\x,T)\le M$. Since $T$ is regular, this implies that
$d(\x,T(x))\le M<\infty$ and therefore $T(x)\ne \emptyset$. Thus $x\in D(T)$. Our assertion follows now from
Theorem  24.1 (c) in [8].

\medskip\noindent
For any $x\in X$ and $\lambda > 0$ consider the following convex lower  semicontinuous function:
$g_{\lambda,x}(z)=\lambda\|z-x\|,\ z\in X$. It is known that
$$\partial g_{\lambda,x}(z)=\{\z\in\X:\|\z\|\le\lambda\ {\rm and}\ \<\z,z-x\>=\lambda\|z-x\|\}\ne\emptyset$$
and in particular
$$\partial g_{\lambda,x}(x)= \lambda \B.$$

\medskip\noindent
{\bf Theorem 2.} A maximal monotone multifunction $T:X\tto \X$ is regular if and only if $T+\partial g_{\lambda,x}$
is maximal monotone for any $x\in X$ and $\lambda > 0$.

\medskip\noindent
{\bf Proof.} The ``only if" part follows from the previous theorem. The ``if" part was essentially proved in [11].
Since it is quite short, we shall repeat it here. To this end, let $(x,\x)\in X\times \X$. We have to show that
$L(x,\x,T)=d(\x,T(x))$. If $L(x,\x,T)=\infty$
there is nothing to prove since, as mentioned earlier, $L(x,\x,T)\le d(\x,T(x))$. So,
assume that $L(x,\x,T)=\lambda < \infty$. A direct computation (see also Lemma 4 in [11]) shows that
$L(x,\x,T+\partial g_{\lambda,x})=0$ which means that $(x,\x)$ is monotonically related
to $T+\partial g_{\lambda,x}$. Since, by hypothesis, $T+\partial g_{\lambda,x}$ is maximal, it follows that
$\x\in (T+\partial g_{\lambda,x})(x)=T(x)+\lambda\B$ and therefore $d(\x,T(x))\le \lambda=L(x,\x,T)$. This completes the proof.

\medskip\noindent
{\bf Corollary 3.} If $T:X\tto \X$ is a maximal monotone multifunction and $D_T$ is either closed and convex
or has non-empty interior then $T$ is regular.

\medskip\noindent
{\bf Proof.} If $D_T$ is closed and convex this follows immediately from Theorem 2 and Voisei's result (the sum
of two maximal monotone multifunctions with closed convex domains that satisfy the usual constraint qualification
is maximal, see [14]). If $D_T$ has non-empty interior our assertion follows from our Theorem 2 and Theorem 9 (i) in [1].

\medskip\noindent
{\bf Strongly maximal monotone multifunctions}

\medskip\noindent
We begin this section by recalling the following definition, due to Simons (see [7], [8]):
a monotone multifunction $T:X\tto\X$ is called {\it strongly maximal} if the following
two conditions are satisfied

(SM1) ``whenever $C\subset X$ is nonempty, convex and $w(X,\X)$-compact, $\x_0\in\X$,
and $(y,\y)\in G(T)$ there exists $x=x_{y,\y}\in C$ such that $\<\y-\x_0,y-x\>\ge 0$"
then there exists $x_0\in C$ such that $\x_0\in T(x_0)$;

(SM2) ``whenever $C\subset \X$ is nonempty, convex and $w(\X,X)$-compact, $x_0\in X$,
and $(y,\y)\in G(T)$ there exists $\x=\x_{y,\y}\in C$ such that $\<\y-\x,y-x_0\>\ge 0$"
then there exists $\x_0\in C$ such that $\x_0\in T(x_0)$.

\medskip\noindent
{\bf Theorem 4.} A maximal monotone multifunction $T:X\tto\X$ is regular if and only if it satisfies
condition (SM2).

\medskip\noindent
{\bf Proof.} It was proved in [11] (Proposition 1) that if $T$ satisfies condition (SM2) then it is regular.
Conversely, assumes that $T$ is regular. We shall adapt a proof of Simons of this assertion
in the case when $T$ is a subdifferential (see for example [8]). Let $x_0\in X$ and $C\subset \X$ be nonempty, convex
and $w(\X,X)$-compact that satisfy the assumption in (SM2). Define a convex, lower
semicontinuous function $f:X\to R$ by $f(x)=\max\<x_0-x,C\>=\max\<x-x_0,-C\>$. It is known
and not difficult to see that

\smallskip
\line{(*)\hfill $\U\in\partial f(x)\ \ {\rm if\ and\ only\ if}\ \ \U\in -C\ {\rm and}\ \<\U,x-x_0\>=f(x).$\hfill}

\smallskip\noindent
On the other hand, if $\y\in T(y)$ and $\x=\x_{y,\y}\in C$ is as in the assumption of (SM2),
then $\<\y-\x,y-x_0\>\ge 0$. It follows that

\smallskip
\line{(**)\hfill $\<\y,x_0-y\>\le \<\x,x_0-y\>\le f(y).$\hfill}

\smallskip\noindent
From (*) and (**) we get that for any $\y\in T(y)$ and $\U\in\partial f(y)$ we have
$$\<\y+\U,y-x_0\>=\<\y,y-x_0\>-\<\U,x_0-y\>=\<\y,y-x_0\>+f(y)\ge 0$$
which means that the pair $(x_0,0)$ is monotonically related to $T+\partial f$. Since $T$ is regular and $\partial f$
is bounded, the maximality of $T+\partial f$ implies that $0=\x_0+\U $ with $\x_0\in T(x_0)$ and
$\U\in\partial f(x_0)=-C$. Thus $\x_0\in T(x_0)\cap C$ and the theorem is proved.

\medskip\noindent
{\bf Corollary 5.} Any strongly maximal monotone multifunction is regular. In particular, any maximal monotone
multifunction whose graph is convex is regular.

\medskip\noindent
{\bf Proof.} The first assertion is obvious while the second one follows from the previous theorem and Theorem 46.1 in [8].

\medskip\noindent
{\bf Enlargements of regular maximal monotone multifunctions}

\medskip\noindent
Let $T:X\tto \X$ be a monotone multifunction. Recall that an {\it enlargement} of $T$
is a multifunction  $E:[0,\infty)\times X\tto \X$ such that $T(x)\subseteq E(\epsi, x)$
for any $x$ and any $\epsi\ge 0$. An enlargement $E:[0,\infty)\times X\tto \X$ is called a {\it
full enlargement} of $T$ if for any $x\in D_T$ and for any $\epsi > 0$  there exists $\delta=\delta(x,\epsi)>0$
such that $T(x)+\delta\B\subseteq T^\epsi$ (see [2], [3] for a general study of enlargements and for
further references and also [6]).

Basically, there are two types of enlargements that are considered. The first one is
defined as follows:
$$E(\epsi,x)=T^\epsi(x)=\{\x\in\X: \<\x-\y,x-y\>\ge -\epsi\|x-y\|,\ (y,\y)\in G(T)\}.$$

\medskip\noindent
{\bf Lemma 6.} With the above notation: $T(x)+\epsi B^*\subseteq T^\epsi(x)$,
for any $x\in D_T$. In particular, $\{T^\epsi\}_{\epsi\ge 0}$ is a full enlargement of $T$.

\medskip\noindent
{\bf Proof.} Let $(x,\x)\in G(T)$ and $\U\in\epsi \B$. Then, for any $(y,\y)\in G(T)$
we have
$$\<\x+\U-\y,x-y\>=\<\x-\y,x-y\>+\<\U,x-y\>\ge-\epsi\|x-y\|$$
which shows that $\x+\U\in T^\epsi(x)$.

\medskip\noindent
{\bf Theorem 7.} Let $T:X\tto\X$ be a regular maximal monotone multifunction, $x\in X$,
$\epsi\ge 0$, and $\x\in T^\epsi(x)$. Then

 \itemitem{(a)} $x\in D_T.$

 \itemitem{(b)} $T^\epsi(x) = T(x)+\epsi\B.$

 \itemitem{(c)} The enlargement $\{T^\epsi\}_{\epsi\ge 0}$ satisfies the following
 monotonicity property:
   $$\<\x-\y,x-y\>\ge-(\epsi+\delta)\|x-y\|,\ \ \x\in T^\epsi(x),\ \y\in T^\delta(y).$$

 \itemitem{(d)} The enlargement $\{T^\epsi\}_{\epsi\ge 0}$ is maximal in the
 following sense: if $(x,\x)\in X\times\X$ and
   $$\<\x-\y,x-y\>\ge-(\epsi+\delta)\|x-y\|\ {\rm for\ any}\ \delta\ge 0 {\rm\ and\ } \y\in T^\delta(y)$$
 then $\x \in T^\epsi(x)$.

 \medskip\noindent
 {\bf Proof.} It is easily seen that
 $$\x\in T^\epsi(x) {\rm\ \ if\ and\ only\ if\ \ } L(x,\x,T)\le\epsi$$
 Thus, since $T$ is regular, if $\x\in T^\epsi(x)$ it follows that $d(\x,T(x))= L(x,\x,T)\le\epsi$ and
 therefore $x\in D_T$ and $\x\in T(x)+\epsi\B$ (since $T(x)$ is weak$^*$ closed and $\epsi \B$ is weak$^*$
 compact). This proves (a) and in view of Lemma 6, also (b).

 To prove (c), let $\x\in T^\epsi(x)$ and $\y\in T^\delta(y)$. In view of part (b),
 $\x=\x_1 + \U$ and $\y=\y_1+\v$ with $\x_1\in T(x),\ \y_1\in T(y),\ \U\in\epsi\B,\ \v\in\delta\B$.
 Then $\|\U-\v\|\le\epsi+\delta$ and
 $$\<\x-\y,x-y\>=\<\x_1+\U-\y_1-\v,x-y\>=\<\x_1-\y_1,x-y\> + \<\U-\v,x-y\>\ge-(\epsi+\delta)\|x-y\|$$
which proves the assertion. Finally, (d) follows from the definition of $T^\epsi(x)$.

\medskip\noindent
{\bf Remark.} One can generalize the definition of $L(x,\x,T)$ as follows:
$$L(x,\x,T^\epsi)=0\vee\sup\Big\{{\<\x-\y,y-x\>\over \|x-y\|}-\epsi-\delta:\,\delta\ge 0, y\ne x, \y\in T^\delta(y)\Big\}\,\ \ \epsi\ge 0$$
It is not difficult to see that if $T$ is regular then $L(x,\x,T^\epsi)=d(\x,T^\epsi(x))$. Indeed,

\smallskip
$\eqalign{L(x,\x,T^\epsi)&=0\vee\sup\Big\{{\<\x-\y,y-x\>\over \|x-y\|}-\epsi-\delta:\,\delta\ge 0, y\ne x, \y\in T^\delta(y)\Big\}\cr
                         &=0\vee\sup\Big\{{\<\x-\z-\U,y-x\>\over\|x-y\|}-\epsi-\delta:\,\delta\ge 0, y\ne x, \z\in T(y),\ \U\in\delta\B\Big\}\cr
						 &=0\vee\sup\Big\{{\<\x-\z,y-x\>\over\|x-y\|}-\epsi+{\<\U,x-y\>\over \|x-y\|}-\delta:\,\delta\ge 0, y\ne x, \z\in T(y),\ \U\in\delta\B\Big\}\cr
						 &=0\vee\sup\Big\{{\<\x-\z,y-x\>\over\|x-y\|}-\epsi:\, y\ne x, \z\in T(y)\Big\}\cr
						 &=0\vee (L(x,\x,T)-\epsi)=0\vee d(\x,T(x))-\epsi=d(\x,T(x)+\epsi \B)=d(\x,T^\epsi(x))}.$

\medskip\noindent
We shall now give a characterization of regularity in terms of the enlargement $\{T^\epsi\}_{\epsi\ge 0}$.

\medskip\noindent
{\bf Theorem 8.} A maximal monotone multifunction $T:X\tto \X$ is regular if and only if $D_{T^\epsi}=D_T$
for any $\epsi>0$.

\medskip\noindent
{\bf Proof.} The ``only if'' part follows from Theorem 7 (a). The ``if'' part will follow from Theorem 2
once we show that $T+\partial g_{\lambda,x}$ is maximal monotone for any $x\in X$ and $\lambda > 0$. To
this end, let $x\in X$ and $\lambda > 0$ and assume that $(z,\z)\in X\times \X$ is monotonically related
to $T+\partial g_{\lambda,x}$. Let $(y,\y)\in G(T)$ and $\U\in \partial g_{\lambda,x}(y)$. We have

\smallskip
\line{\qquad $\<\z-\y,z-y\>= \<\z-\y-\U,z-y\>+\<\U,z-y\>\ge 0-\lambda\|z-y\|=-\lambda\|z-y\|$\hfill}

\smallskip\noindent
which shows that $\z\in T^\lambda(z)$ and in particular $z\in D_{T^\lambda}=D_T$. Theorem 24.1 (c) in [8]
proves that $T+\partial g_{\lambda,x}$ is maximal monotone and this finishes the proof.

\medskip\noindent
We shall now turn our attention to another type of enlargements which were studied in
[3] and [6]. If $T:X\tto \X$ is a monotone multifunction and $x\in X$ define
$$E(\epsi,x)=T_\epsi(x)=\{\x\in\X:\<\x-\y,x-y\>\ge-\epsi\ {\rm for \ any}\ (y,\y)\in G(T)\}.$$
It is worth mentioning that this enlargement belongs to the class
I\hskip-1.6pt E$(T)$ introduced in [10] while the enlargement $\{T^\varepsilon\}$ considered earlier does not. 

\medskip\noindent
{\bf Theorem 9.} If $T:X\tto \X$ is a regular maximal monotone multifunction, $\epsi > 0$ and
$T_\epsi(x)\ne\emptyset$ then  $x\in \overline{D_T}$, that is $D_{T_\epsi}\subseteq \overline{D_T}$.

\medskip\noindent
 {\bf Proof.} Assume not. Then there exists $\delta>0$ such that $\|x-y\|>\delta$ for
 any $y\in D_T$. Let $\x\in T_\epsi(x)$. Then
$$\<\x-\y,x-y\>\ge-\epsi\ge -{\epsi\over\delta}\|x-y\|\ {\rm for \ any}\ (y,\y)\in G(T)$$
and therefore $\x\in T^{\epsi/\delta}(x)$. By Theorem 7 (a), $x\in D_T$, which is a
contradiction. It follows that $x\in \overline{D_T}$.

\medskip\noindent
{\bf Remark.} A particular case of this result (when $X$ is reflexive) was proved in [3].

\bigskip\noindent
{\bf Acknowledgement.} We would like to thank the referee for carefully reading the manuscript and for making several useful comments.

\bigskip\noindent
{\bf References}

\medskip
\item{[1]} J. M. Borwein: Maximality of sums of two maximal monotone operators in general Banach spaces, Proc.\ Amer.\ Math.\ Soc.,
              135 (2007) 3917--3924.

\item{[2]} R. S. Burachik and A. N. Iusem: On non-enlargeable and fully enlargeable monotone operators, J. of Convex Analysis,
              13 (2006) 603–-622.

\item{[3]} R. S. Burachik and B. F. Svaiter: $\epsi$-enlargements of maximal monotone operators in Banach spaces, Set-Valued Analysis,
              7 (1999) 117–-132.

\item{[4]} M. Marques Alves and B. F. Svaiter: A new old class of maximal monotone operators, J. of Convex Analysis,
              16 (4) (2009), to appear.

\item{[5]} R. R. Phelps: Lectures on Maximal Monotone Operators, Extracta Mathematicae 15 (1997), 193--230.

\item{[6]} J. P. Revalski and M. Th\'era: Enlargements and sums of monotone operators, Nonlinear Ana. 48 (2002) 505--519.

\item{[7]} S. Simons: Minimax and Monotonicity, Lecture Notes in Mathematics 1693 (1998), Springer-Verlag.

\item{[8]} S. Simons: From Hahn-Banach to Monotonicity, Second Edition, Lecture Notes in Mathematics 1693 (2008), Springer-Verlag.

\item{[9]} S. Simons: The least slope of a convex function and the maximal monotonicity of its subdifferential, J. Optim.\ Theory 71 (1991) 127--136.

\item{[10]} B. F. Svaiter: A family of enlargements of maximal monotone operators, Set-Valued Analysis,
              8 (2000) 311–-328.

\item{[11]} A. Verona and M. E. Verona: Regular maximal monotone operators, Set-Valued Anal.\ 6 (1998), 302--312.

\item{[12]} A. Verona and M. E. Verona: Regular maximal monotone operators and the sum theorem, J. of Convex Analysis,
          7 (2000) 115--128.

\item{[13]} A. Verona and M. E. Verona: Regularity and the Br\o ndsted-Rockafellar properties of maximal monotone operators,
           Set-Valued Anal.\ 14 (2006), 149--157.

\item{[14]} M. D. Voisei: A maximality theorem for the sum of maximal monotone operators in non-reflexive Banach spaces,
          Math.\ Sci.\ Res.\ J. 10 (2006) 36--41.

\item{[15]} M. D. Voisei and C. Z\u alinescu: Strongly-representable monotone operators, Preprint (2008).
\bye